\documentclass[12pt]{amsart}

\usepackage{amsmath, amscd, graphicx, latexsym, hyperref, rlepsf, times, pdfsync}

\newtheorem{Thm}{Theorem}[section]
\newtheorem{Lem}[Thm]{Lemma} 
\newtheorem{Prop}[Thm]{Proposition}

\theoremstyle{definition}
\newtheorem{Def}[Thm]{Definition}
\newtheorem{exa}[Thm]{Example}

\theoremstyle{remark}
\newtheorem*{Rem}{Remark}
\newtheorem*{Ack}{Acknowledgements}

\numberwithin{equation}{section}

\newcommand{\pob}{partial open book decomposition }

\def\p{\partial}

\def\v{\vskip.12in}
\def\a{\alpha}
\def\A{$\boldmath$\alpha$\unboldmath$}
\def\B{$\boldmath$\beta$\unboldmath$}
\def\Atek{\boldmath$\alpha \ $\unboldmath}
\def\Btek{\boldmath$\beta \ $\unboldmath}
\def\b{\beta}
\def\d{\delta}
\def\G{\Gamma}
\def\g{c}

\def\S{\Sigma}

\begin{document}

\title[On the relative Giroux correspondence]{On the relative Giroux correspondence}

\author{Tolga Etg\"u}

\author{Burak Ozbagci}

\begin{abstract}
Recently, Honda, Kazez and Mati\'{c} described an 
adapted \pob of a compact contact $3$-manifold with convex
boundary by generalizing the work of Giroux in the closed case. 
They also implicitly established  a one-to-one
correspondence between isomorphism classes of partial open book
decompositions modulo positive stabilization and isomorphism
classes of compact contact $3$-manifolds with convex boundary. 
In this expository article we explicate  the relative version of Giroux correspondence. 
\end{abstract}

\address{Department of Mathematics \\ Ko\c{c} University \\ Istanbul, Turkey}
\email{tetgu@ku.edu.tr} \email{bozbagci@ku.edu.tr}
\subjclass[2000]{}


\thanks{The first author  was
partially supported by a GEBIP grant of the Turkish Academy of
Sciences and a CAREER grant of the Scientific and Technological
Research Council of Turkey.}
\thanks{The second author was partially supported by the
research grant 107T053 of the Scientific and Technological
Research Council of Turkey and the
Marie Curie International Outgoing Fellowship 236639.}

\v \v \v

\maketitle

\setcounter{section}{-1}


\section{Introduction}

Let $(M,\G)$ be a balanced sutured $3$-manifold and let $\xi$ be a
contact structure on $M$ with convex boundary whose dividing set
on $\p M$ is isotopic to $\G$. Recently, Honda, Kazez and
Mati\'{c} \cite{hkm1}  introduced an invariant of the contact
structure $\xi$ which lives in the sutured Floer homology group
defined by Juh\'{a}sz  \cite{juh}. This invariant is a relative
version of the contact class in Heegaard Floer homology in the
closed case as defined by Ozsv\'{a}th and Szab\'{o} \cite{os} and
reformulated in \cite{hkm}. Both the original definition in
\cite{os} and the reformulation of the contact class by Honda,
Kazez and Mati\'{c} are based on the so called Giroux
correspondence \cite{g} which is a one-to-one correspondence
between open book decompositions modulo
positive stabilization and isotopy classes of contact
structures on closed $3$-manifolds.

In order to adapt their reformulation \cite{hkm} of the contact
class to the case of a contact manifold $(M, \xi)$  with convex
boundary, Honda, Kazez and Mati\'{c} described in \cite{hkm1}, a
\pob of $M$ (adapted to $\xi$)  by generalizing the work of Giroux
in the closed case. This description coupled with Theorem 1.2 (and
the subsequent discussion) in \cite{hkm1} induces a map from
isomorphism classes of compact contact $3$-manifolds with convex
boundary to isomorphism classes of partial open book
decompositions modulo positive stabilization. Here we spell out the inverse of this map, 
by describing a compact contact
$3$-manifold with convex boundary compatible with an
\emph{abstract} partial open book decomposition. To define a contact
structure compatible with an abstract partial open book
decomposition we chose to mimic the analogous result of
Torisu \cite{to} (rather than adapting the construction of Thurston and
Winkelnkemper \cite{tw})   which conveniently allowed us to keep track of
the dividing set on the boundary. Consequently,  one obtains a relative version
of Giroux correspondence which is due to Honda, Kazez and Mati\'{c}.

\begin{Thm}\label{giroux} There is a one-to-one
correspondence between isomorphism classes of partial open book
decompositions modulo positive stabilization and isomorphism
classes of compact contact $3$-manifolds with convex boundary.
\end{Thm}

The relative Giroux correspondence helps understand the geometric properties of contact $3$-manifolds using partial open books, e.g. if the monodromy of a corresponding partial open book is not right-veering, then the contact structure is overtwisted. It also plays a critical role in the definition of the (relative) contact invariant in sutured Floer homology which helps to analyze the contact invariant of a closed manifolds in terms of the relative contact invariants of certain compact pieces. In \cite{ghv}, it is proved that the contact invariant vanishes in the presence of Giroux torsion using some properties of the relative invariant. 

The paper is organized as follows: In Section~\ref{def} we give the definition of an abstract \pob $(S,P,h)$, construct a
balanced sutured manifold $(M,\G)$ associated to $(S,P,h)$, and construct a (unique) \emph{compatible} contact structure $\xi$ on
$M$ which makes $\p M$ convex with a dividing set isotopic to
$\G$. In Section~\ref{relative} we prove Theorem~\ref{giroux} after reviewing
the related results due to Honda, Kazez and Mati\'{c} \cite{hkm1}. In the last section we provide examples of abstract partial open books compatible with some basic contact $3$-manifolds with boundary. 

 The reader is advised to turn to Etnyre's
notes \cite{e} for the related material on contact topology of
$3$-manifolds.

\begin{Ack}
We would like to thank Andr\'as 
Stipsicz, Sergey Finashin and John Etnyre for valuable comments on a draft of
this paper. We also thank the anonymous referee for helpful remarks, especially for the remark in the footnote for Proposition~\ref{torisu}. 
\end{Ack}

\section{Partial open books, sutured manifolds  and contact structures}\label{def}

\begin{Def} \label{Pob} An abstract  \pob is a triple $(S,P,h)$ satisfying the
following conditions:

$(1)$ $S$ is a compact oriented connected surface with $\p S \neq
\emptyset$,

$(2)$ $P= P_1 \cup P_2 \cup \ldots \cup P_r$ is a proper (not
necessarily connected) subsurface of $S$ such that $S$ is obtained
from $\overline{S \setminus P}$ by successively attaching $1$-handles $P_1,
P_2, \ldots, P_r$,

$(3)$ $h:P \to S$ is an embedding such that $h|_A =$ identity, where
$A=\p P\cap \p S$.

\end{Def}

\begin{figure}[ht]

  \relabelbox \small {
  \centerline{\epsfbox{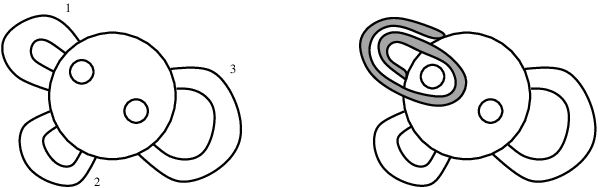}}}
  \relabel{1}{{$P_1$}}
  \relabel{2}{{$P_2$}}
 \relabel{3}{{$P_3$}}

  \endrelabelbox
        \caption{An example of $S$ and $P$ satisfying the conditions in Definition~\ref{Pob}:
        $\overline{S \setminus P}$ is a twice punctured disk, $r=3$, and $h$ is the embedding
        which is identity on $P_2$ and $P_3$, and the image of $P_1$ is the shaded region indicated in the figure on the right.}
        \label{def1}
\end{figure}

\begin{figure}[ht]
  \relabelbox \small {
  \centerline{\epsfbox{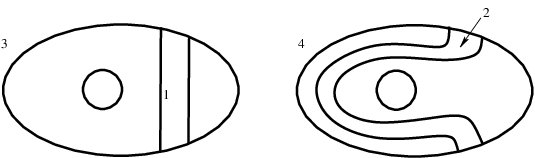}}}

  \relabel{2}{$h(P)$}

   \relabel{4}{$S$}

   \relabel{1}{$P$}
   \relabel{3}{$S$}
\endrelabelbox

        \caption{Another example of an abstract partial open book.}
        \label{ex4-1}
\end{figure}

\begin{figure}[ht]
  \relabelbox \small {
  \centerline{\epsfbox{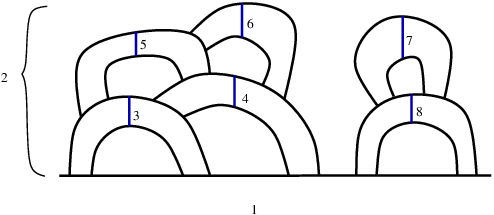}}}

  \relabel{1}{$S\setminus P$}
  \relabel{2}{$P$}
  \relabel{3}{$a_1$}
  \relabel{4}{$a_2$}
  \relabel{5}{$a_3$}
 \relabel{6}{$a_4$}
  \relabel{7}{$a_5$}
  \relabel{8}{$a_6$}

\endrelabelbox

        \caption{A basis of $H_1(P,A)$: cocores $a_1 ,a_2 ,\ldots,a_6 $ of the
        $1$-handles in $P$}
        \label{cocores}
\end{figure}

\begin{Rem}

Figures~\ref{def1} and \ref{ex4-1} present simple examples of partial open book decompositions. 
It follows from the definition that $A$ is a $1$-manifold
with nonempty boundary (but it may have closed components as in Figure~\ref{cocores}) and $\overline{\p P \setminus A}$ is a
nonempty set consisting of some arcs (but no closed components).
The connectedness condition on $S$ is not essential, but
simplifies the discussion.

\end{Rem}

\v

We now  briefly turn our attention to sutured manifolds which was introduced by Gabai \cite{ga} to study foliations. 
 A sutured manifold $(M,\G)$ is a compact oriented $3$-manifold
with nonempty boundary, together with a compact subsurface $\G=
A(\G) \cup T(\G) \subset \p M$, where $A(\G)$ is a union of
pairwise disjoint annuli and $T(\G)$ is a union of tori. Moreover
each component of $\p M \setminus \G$ is oriented, subject to the
condition that  whether or not the orientation agrees with the orientation 
induced as the boundary of $M$ changes  every time we nontrivially
cross $A(\G)$.  Let $R_+ (\G)$ (resp. $R_- (\G)$) be the open
subsurface of $\p M \setminus \G$ on which the orientation agrees
with (resp. is the opposite of ) the boundary orientation on $\p
M$. A sutured manifold $(M,\G)$ is balanced if $M$ has no closed
components, $\pi_0(A(\G)) \to \pi_0(\p M) $ is surjective, and
$\chi(R_+ (\G))= \chi (R_- (\G))$ on every component of $M$. It
turns out that if $(M,\G)$ is balanced, then $\G=A(\G)$ and every
component of $\p M$ nontrivially intersects $\G$. Since all the
sutured manifolds that we will deal with in this paper are
balanced, we will think of $\G$ as a set of \emph{oriented curves}
on $\p M$ by identifying each annulus in $\G$ with its core
circle. Here we orient $\Gamma$ as the boundary of $R_+ (\G)$.

We now emphasize the relation between dividing sets and sutures.
Let $\xi$ be a contact structure on a compact oriented
$3$-manifold $M$ whose dividing set on the convex boundary $\p M$
is denoted by $\G$. Then it is fairly easy to see that $(M,
\Gamma)$ is a \emph{balanced} sutured manifold (with annular
sutures) via the identification we mentioned above. Conversely,
given a balanced sutured manifold $(M,\G)$, there exists a contact
structure $\xi$ on $M$ which makes $\p M$ convex and realizes $\G$
as its diving set on $\p M$. However one should keep in mind
that the contact structure is not uniquely determined and cannot
always be chosen to be tight.

Given a \pob $(S,P,h)$, we construct a sutured manifold $(M,\G)$
as follows: Let $$H=(S \times [-1,0])/\sim
$$ where $(x,t) \sim (x,t')$ for $x \in \p S$ and $t, t' \in
[-1,0]$. It is easy to see that $H$ is a solid handlebody whose
oriented  boundary is the surface $S \times \{0 \} \cup - S \times
\{-1\}$ (modulo the relation $(x,0) \sim (x,-1)$ for every $x \in
\p S$). Similarly let $$N=(P \times [0,1])/\sim$$ where $(x,t)
\sim (x,t')$ for $x \in A$ and $t, t' \in [0,1]$. Since $P$ is not
necessarily connected $N$ is not necessarily connected. Observe
that each component of $N$ is also a solid handlebody. The
oriented boundary of $N$ can be described as follows: Let the arcs
$\g_1, \g_2, \ldots, \g_n$ denote the connected components of
$\overline{\p P \setminus  A}$. Then, for $1 \leq i \leq n$, the
disk $D_i= (\g_i \times [0,1])/\sim$ belongs to $\p N$. Thus part
of $\p N$ is given by the disjoint union of $D_i$'s. The rest of
$\p N$ is the surface $P \times \{1 \} \cup - P \times \{0\}$
(modulo the relation $(x,0) \sim (x,1)$ for every $x \in A$).

\begin{figure}[ht]
  \relabelbox \small {
  \centerline{\epsfbox{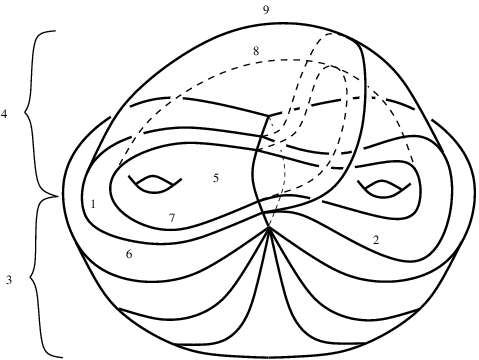}}}

  \relabel{1}{$P$}
  \relabel{2}{$h(P)$}
  \relabel{3}{$H$}
  \relabel{4}{$N$}
  \relabel{5}{$S$}
 \relabel{6}{$c_1$}
  \relabel{7}{$c_2$}
  \relabel{8}{$D_2$}
  \relabel{9}{$D_1$}

\endrelabelbox

        \caption{A partial open book decomposition: $M$ as the union of $N$ and $H$}
        \label{pobook}
\end{figure}

Let $M= N \cup H$ where we glue these manifolds by identifying
$P\times \{0\} \subset \p N $ with $P \times \{0\} \subset \p H $
and $P \times \{1\} \subset \p N $ with $h(P) \times \{-1\}
\subset \p H $. Since the gluing identification is orientation
reversing $M$ is a compact oriented $3$-manifold with oriented boundary
$$\p M =
(S\setminus P ) \times \{ 0\}
\cup -(S \setminus h(P)) \times \{ -1\}
\cup (\overline{\p P \setminus  A}) \times [0,1]
$$
(modulo the identifications given above).

\begin{Def} \label{part} If a compact $3$-manifold $M$ with boundary is
obtained from $(S,P,h)$ as discussed above, then we call the triple
$(S,P,h)$ a partial open book decomposition of $M$.

\end{Def}

We define the suture $\G$ on $\p M$ as the set of closed curves
(see Remark~\ref{curve}) obtained by gluing the arcs $\g_i \times
\{1/2 \} \subset \p N$, for $1 \leq i \leq n$, with the arcs in
$(\overline{\p S \setminus \p P}) \times \{0\} \subset \p H$,
hence as an oriented simple closed curve and modulo
identifications
$$ \G =
(\overline{\p S \setminus \p P}) \times \{0\}
\cup
- (\overline{\p P \setminus A}) \times \{1/2\}
\ .$$
\begin{Rem} \label{curve} If a sutured manifold $(M,\G)$ has only
annular sutures, then it is convenient to refer to the set of core
circles of these annuli as $\G$.
\end{Rem}

\begin{Def} \label{abs} The sutured manifold $(M,\G)$ obtained
from  a \pob $(S,P,h)$ as described above is called the sutured
manifold associated to $(S,P,h)$.
\end{Def}

\begin{Def} [\cite{juh}] A sutured manifold $(M,\G)$ is balanced if $M$ has no
closed components, $\pi_0(A(\G)) \to \pi_0(\p M) $ is surjective,
and $\chi(R_+ (\G))= \chi (R_- (\G))$ on every component of $M$.
\end{Def}

\begin{Rem} \label{annuli} It follows that if $(M,\G)$ is balanced,
then $\G=A(\G)$ and every component of $\p M$ nontrivially
intersects the suture $\G$.
\end{Rem}

{\Lem The sutured manifold $(M,\G)$
associated to a \pob $(S,P,h)$ is balanced.
}

\begin{proof} It is clear that $M$ is connected since we assumed that $S$ is
connected. We observe that $\p M \neq \emptyset$ since $P$ is a
proper subset of $S$ by our definition. In fact, $\p M$ can be described starting from the connected surface $\p H = S \times \{ 0\} \cup -S \times \{ -1\}$: Let $\kappa_j$ be $a_j \cup h(a_j)$, where $a_j$ is the cocore of the 1-handle $P_j$ in $P$ (see Figure~\ref{cocores} for suitable $a_j$'s). Then $\p M$ is obtained by cutting $\p H$ along $\kappa_j$'s and capping off each resulting boundary by a disk $D_i= (\g_i \times [0,1])/\sim$ \! \! for some $i$. From this description it is clear that  every component of $\p M$
contains a $\g_i \times
\{1/2 \} \subset \G$ and therefore $\pi_0(A(\G)) \to \pi_0(\p M) $
is surjective. Now let $R_+ (\G)$ be the open subsurface in $\p M$
obtained by gluing
$$((S \setminus \p S) \setminus P) \times \{0\} \subset \p H\;\;  \mbox{and} \;\;
\cup_{i=1}^{n}\; (\g_i \times [0,1/2))/\sim \;\; \subset \p N$$
and $R_- (\G)$ be the open subsurface in $\p M$ obtained by gluing
$$((S \setminus \p S)\setminus h(P)) \times \{-1\} \subset \p H\;\;
\mbox{and}\;\; \cup_{i=1}^{n}\; (\g_i \times (1/2,1])/\sim \;\;
\subset \p N$$ under the gluing map that is used to construct $M$.
Since $h : P \to S $ is an embedding we have $\chi(P)= \chi
(h(P))$ and it follows that $\chi(R_+ (\G))= \chi (R_- (\G))$.
\end{proof}

The following result is inspired by Torisu's work \cite{to} in the
closed case.

{\Prop \label{torisu} Let $(M,\G)$ be the balanced sutured
manifold associated to a \pob $(S,P,h)$. Then  there exists a
contact structure $\xi$ on $M$ satisfying the following
conditions:

$(1)$ $\xi$ is tight when restricted to $H$ and $N$,

$(2)$ $\p H$ is a convex surface in $(M, \xi)$ whose dividing set
is $\p S \times \{0 \}$,

$(3)$ $\p N $ is a convex surface in $(M, \xi)$ whose dividing set
is $\p P \times \{ 1/2 \}$.

Moreover such $\xi$ is unique up to isotopy.}

\begin{proof} We will prove that there is a unique tight contact structure
(up to isotopy) on $H$ and $N$ with the given boundary
conditions, using arguments along similar lines.\footnote{In fact, one can prove a general existence and uniqueness theorem using an explicit contact form $\lambda + dt$ on $\Sigma \times [0,1] / \sim $ \!  , for any surface $\Sigma$ with boundary, where $\lambda$ is a primitive of a volume form on $\Sigma$ that is standard near the boundary. It can be argued that this contact form gives a tight contact structure making the boundary convex with dividing set $\p \Sigma \times \{ 1/2\}$.}\label{referee}
 Once we have these contact structures on $H$ and $N$, since the dividing sets on $\p H$ and $\p N$ agree
on the subsurface along which we glue
$H$ and $N$, we obtain a unique contact
structure (up to isotopy) on $M$ satisfying the above conditions, by gluing together the contact structures on these pieces. 

To prove the existence of tight contact structures on $H$ and $N$ with prescribed  dividing sets we simply consider $H$ and $N$ embedded in the closed contact $3$-manifold $(Y, \xi')$ supported by the open book $(S, id)$ and appeal to the closed case (see \cite{to} and \cite[Lemma 4.4]{e}). For $H$, observe that 
$$ H = ( S \times [-1, 0] ) / \sim \ \ \subset ( S \times [-1, 1] ) / \sim  \ \ = Y  \ , $$
where the equivalence relation $\sim$ is given by, $(x,t) \sim (x, t')$ for $x \in \p S$ and $t,t' \in [-1, 1]$, and $(s,-1) \sim (s, 1)$ for $s \in S$. The contact structure $\xi'$ is Stein fillable by \cite{g}, hence tight by \cite{eg}, and therefore its restriction to $H$ is also tight. In fact, $\p H$ is convex with respect to $\xi'$ with dividing set $\p S \times \{ 0 \}$ (see Lemma 4.4 in \cite{e}). 
Similarly, $N$ trivially embeds in $H$ since $\p P \times \{ 1/2 \}$ is the union of $A
\times \{0\}$ and the arcs $\g_i \times \{1/2 \} $, for $1\leq
i\leq n$. So $\xi'$ restricts to a tight contact structure on $N$. To identify its dividing set we
first observe that the
dividing set on $P \times \{1 \} \cup - P \times \{0\}=
\p N \cap \p H$ is the set $A \times \{ 0\}= \p N \cap (\p S \times \{ 0\})$.
The rest of $\p N$ consists of the
disks $D_i=(\g_i \times [0,1])/\sim$. Each one of these disks
can be made convex so that the dividing set is a single arc
since its boundary intersects the dividing set twice. It follows that the
dividing set on $\p N$ is as required after rounding the edges.

In order to prove the uniqueness for $H$, as in Lemma 4.4 in \cite{e}, we take a set $\{ d_1 
,\ldots,d_p \}$ of properly embedded pairwise disjoint arcs in $S$
whose complement is a single disk. (It follows that the set $\{
d_1 ,d_2 ,\ldots,d_p \}$ represents a basis of $H_1(S,\p S)$.) For
$ 1\leq k\leq p$, let $\d_k$ denote the closed curve on $\p H$
which is obtained by gluing the arc $d_k$ on $S \times \{0\}$ with
the arc $d_k$ on $S \times \{-1\}$.  Then we observe that $\{\d_1,
\d_2, \ldots, \d_p \}$ is a set of homologically linearly
independent closed curves on $\p H$ so that $\d_k$ bounds a
compressing disk $D^\d_k=(d_k \times [0,-1]) / \sim$ in $H$. It is
clear that when we cut $H$ along $D^\d_k$'s (and smooth the
corners) we get a $3$-ball $B^3$. Moreover $\d_k$ intersects the
dividing set twice by our construction. Now we put each $\d_k$
into Legendrian position (by the Legendrian realization principle
\cite{hon}) and make the compressing disk $D^\d_k$ convex
\cite{g1}. The dividing set on $D^\d_k$ will be an arc connecting
two points on $\p D^\d_k=\d_k$. Then we cut along these disks and
round the edges (see \cite{hon}) to get a connected dividing set
on the remaining $B^3$. Consequently,  Theorem~\ref{elia} due to Eliashberg 
(although stated in different terms in \cite{eli}) 
implies the uniqueness of a tight contact
structure on $H$ with the assumed boundary conditions. 
Recall that  a standard contact $3$-ball is a tight contact $3$-ball with convex boundary whose dividing set is connected.

\begin{Thm}  [Eliashberg] \label{elia} There is a unique standard contact $3$-ball. \end{Thm}

The proof of the uniqueness of such a
tight contact structure on $N$ follows a similar line. Instead of a basis of $H_1(S, \p S)$
we take suitable cocores $\{ a_1 , \dots, a_r \}$ of the 1-handles $P_j$'s in $P$ to 
get a basis of $H_1(P,A)$  (see Figure~\ref{cocores} for an
example). Then one can proceed as in the proof  given above for the handlebody $H$.
\end{proof}

Proposition~\ref{torisu} leads to the following definition of
compatibility of a contact structure and a partial open book
decomposition.

\begin{Def} \label{pOb} Let $(M,\G)$ be the balanced sutured manifold
associated to a \pob $(S,P,h)$. A contact structure $\xi$ on
$(M,\G)$ is said to be compatible with $(S,P,h)$ if it is isotopic to a contact structure satisfying
conditions $(1),(2)$ and $(3)$ stated in
Proposition~\ref{torisu}.
\end{Def}

\begin{Def} \label{iso} Two partial open book decompositions $(S,P,h)$
and $(\widetilde{S}, \widetilde{P}, \widetilde{h})$ are isomorphic
if there is a diffeomorphism $f: S \to \widetilde{S}$ such that
$f(P) = \widetilde{P}$ and $\widetilde{h} = f \circ h \circ
(f^{-1})|_{\widetilde{P}}$.
\end{Def}

\begin{Rem} \label{unique} It follows from Proposition~\ref{torisu} that
every \pob has a unique compatible contact structure, up to
isotopy, on the balanced suture manifold associated to it, such
that the dividing set of the convex boundary is isotopic to the
suture. Moreover if $(S,P,h)$ and $(\widetilde{S}, \widetilde{P},
\widetilde{h})$ are isomorphic partial open book decompositions,
then the associated compatible contact $3$-manifolds $(M, \G,
\xi)$ and $(\widetilde{M}, \widetilde{\G}, \widetilde{\xi})$ are also
isomorphic.

\end{Rem}

\begin{Def} Let $(S,P,h)$ be a partial open book decomposition. A partial open book
decomposition $(S',P',h')$ is called a positive stabilization of
$(S,P,h)$ if there is a properly embedded arc $s$ in $S$ such that
\begin{itemize}
\item $S'$ is obtained by attaching a $1$-handle to $S$ along $\p
s$, \item $P'$ is defined as the union of $P$ and the attached
$1$-handle, \item $h' = R_{\sigma} \circ h$, where the extension
of $h$ to $P'$ by the identity is also denoted by $h$, and $R_{\sigma}
$ denotes the right-handed Dehn twist along the closed curve
$\sigma$ which is the union of $s$ and the core of the
attached $1$-handle.
\end{itemize}
\end{Def}

The effect of positively stabilizing a \pob on the associated
sutured manifold and the compatible contact structure is taking a
connected sum with $(S^3, \xi_{std})$ away from the boundary. We will prove this statement in Lemma~\ref{invariance} 
and the notion of sutured Heegaard diagram will be helpful in our argument.  So we digress to review basic definitions and properties of
Heegaard diagrams of sutured manifolds (cf. \cite{juh}).  

A sutured Heegaard diagram is given by $(\S, \A, \B)$, where the
Heegaard surface $\S$ is a compact oriented surface with nonempty
boundary and
 $\A= \{\a_1, \a_2, \ldots, \a_m \} $ and $\B = \{\b_1, \b_2, \ldots, \b_n \}
$ are two sets of pairwise disjoint simple closed curves in $\S
\setminus \p \S$. Every sutured Heegaard diagram  $(\S, \A, \B)$,
uniquely defines a sutured manifold $(M, \G)$ as follows: Let $M$
be the $3$-manifold obtained from $\S \times [0,1]$ by attaching
$3$-dimensional $2$-handles along the curves $\a_i \times \{0\}$
and $\b_j \times \{1\}$ for $i=1, \ldots, m$ and $ j= 1, \ldots,
n$. The suture $\G$ on $\p M$ is defined by the set of curves $\p
\S \times \{ 1/2 \}$ (see Remark~\ref{curve}).

In \cite{juh}, Juh\'asz proved that if  $(M, \G)$ is defined by
$(\S, \A, \B)$, then $(M, \G)$ is balanced if and only if
$|\A|=|\B|$, the surface $\S$ has no closed components and both
\Atek and \Btek consist of curves linearly independent in $H_1 (
\S , \mathbb{Q})$. Hence a sutured Heegaard diagram $(\S, \A, \B)$
is called balanced if it satisfies the conditions listed above. We
will abbreviate balanced sutured Heegaard diagram as balanced
diagram.

A partial open book decomposition of $(M,\G)$ gives a sutured
Heegaard diagram $(\S,\A,\B)$ of $(M,-\G)$ as follows: Let
$$\Sigma=P \times \{ 0\}
 \cup - S \times \{-1\} / \sim \ \subset \p H$$
be the Heegaard surface. Observe that, modulo identifications,
$$\p \S =
(\overline{\p P \setminus A}) \times \{0\} \cup -(\overline{\p
S \setminus \p P}) \times \{-1\} \simeq -\G \ . $$ As in the proof
of Proposition~\ref{torisu}, let $a_1, a_2, \ldots , a_r$ be
properly embedded pairwise disjoint arcs in $P$ with endpoints on
$A$ such that $S \setminus \cup_j a_j$ deformation retracts onto
$\overline{S\setminus P}$. Then define two families $\A = \{
\alpha_1, \a_2, \ldots , \alpha_r \}$ and $\B = \{ \beta_1, \b_2,
\ldots , \beta_r\}$ of simple closed curves in the Heegaard
surface $\S$ by $$\alpha_j = a_j \times \{ 0\} \cup a_j \times
\{-1\} / \sim \ \mbox{ and } \ \beta_j = a_j \times \{ 0\} \cup h(a_j)
\times \{ -1 \} / \sim \ \ .$$
$(\Sigma, \A, \B)$ is a sutured Heegaard diagram of $(M,-\G)$.
Here the suture is $-\G$ since $\p \S$ is isotopic to $-\G$.

\begin{Lem}\label{invariance}
The balanced sutured manifold associated to a partial open book decomposition and
the compatible contact structure are invariant under positive stabilization.
\end{Lem}
\begin{proof}
Let $(S,P,h)$ be a partial open book decomposition of $(M,\G)$,
$s$ be a properly embedded arc in $S$, and $(S',P',h')$ be the
corresponding positive stabilization of $(S,P,h)$. Consider the
sutured Heegaard diagram $(\S,\A, \B)$ of $(M,-\G)$ given by
$(S,P,h)$ using properly embedded disjoint arcs $a_1, a_2 , \dots , a_r$ in
$P$.

Let $a_0$ be the cocore of the 1-handle attached to $S$ during
stabilization. The endpoints of $a_0$ are on $A'=\p P' \cap \p S'$
and $S' \setminus \cup_{j=0}^r a_j $ deformation retracts onto $S'
\setminus P' = S \setminus P$. Using the properly embedded disjoint arcs
$a_0, a_1,a_2, \dots, a_r$ in $P'$ we get a sutured Heegaard
diagram $(\S',\A',\B')$ of $(M',-\G')$, where $(M',\G')$ is the
sutured manifold associated to $(S',P',h')$. Observe that $\A' =
\{ \a_0\} \cup$ \Atek \!, $\B' = \{ \b_0 \} \cup$ \Btek \!, and
$$\S' = P' \times \{ 0\} \cup -S' \times \{ -1\} / \!\!\sim \ \ \cong T^2 \# \S \ \ . $$
Since $h'$ is a right-handed Dehn twist along $\sigma$ composed
with the extension of $h$ which is identity on $P' \setminus P$,
$\a_0$ is disjoint from
every  $\b_j$ with $j > 0$. Therefore $(\S',\A',\B')$ is a
stabilization of the Heegaard diagram $(\S,\A,\B)$, and
consequently $(M', \G') \cong (M, \G)$. The contact structure
$\xi'$ compatible with $(S',P',h')$ is contactomorphic  to $\xi$
since $\xi'$ is obtained from $\xi$ by taking a connected sum 
with
$(S^3, \xi_{std})$ away from the boundary. This can be seen as in the closed 
case, and holds essentially because of the fact that the abstract open book with an annulus page 
and monodromy given by a right-handed Dehn twist (which is the one that gives the genus-1
Heegaard decomposition with a single \Atek-curve that intersects the single \Btek-curve geometrically once) is compatible with the standard 
contact structure on $S^3$.
\end{proof}

\section{Relative Giroux correspondence}\label{relative}

The following theorem is the key to
obtaining a description of a partial open book decomposition of
$(M,\G, \xi)$ in the sense of Honda, Kazez and Mati\'{c}.

{\Thm [\cite{hkm1}, Theorem 1.1]  \label{honkazmat} Let $(M,\G)$
be a balanced sutured manifold and let $\xi$ be a contact
structure on $M$ with convex boundary whose dividing set $\G_{\p
M}$ on $\p M$ is isotopic to $\G$. Then there exist a Legendrian
graph $K \subset M $ whose endpoints lie on $\G \subset \p M$ and
a regular neighborhood $N(K) \subset M $ of $K$ which satisfy the
following:

\begin{itemize}
\item[](A) \ (i) $T=\overline{\p N(K) \setminus \p M}$ is a convex surface
with Legendrian boundary.
\begin{itemize}
\item[](ii) For each component $\gamma_i$ of
$\p T$, $\gamma_i \cap \G_{\p M}$ has two connected components.
\item[](iii) There is a system of pairwise disjoint compressing disks
$D^\a_j$ for $N(K)$ so that $\p D^\a_j$ is a curve on $T$
intersecting the dividing set $\G_T$ of $T$ at two points and each
component of $N(K) \setminus \cup_j D^\a_j$ is  a standard contact
$3$-ball, after rounding the edges.
\end{itemize}
\item[](B) \ (i) Each component $H$ of $\overline{M \setminus N(K)}$ is a
handlebody (with convex boundary).
\begin{itemize}
\item[](ii) There is a system of
pairwise disjoint compressing disks $D^{\d}_k$ for $H$ so that
each $\p D^{\d}_k$ intersects the dividing set $\G_{\p H}$ of $\p
H$ at two points and $H \setminus \cup_k D^{\d}_k$ is a standard
contact $3$-ball, after rounding the edges.
\end{itemize}
\end{itemize}}

\vspace{0.1in}

Based on Theorem~\ref{honkazmat}, Honda, Kazez and Mati\'{c}
describe a \pob on $(M,\G)$ in \cite[Section 2]{hkm1}. In this paper, for the sake of simplicity and without
loss of generality, we will assume that $M$ is connected. As a
consequence $M\setminus N(K)$ in Theorem~\ref{honkazmat} is also
connected.

We claim that the description in \cite{hkm1} gives
a \pob $(S,P,h)$; that the balanced sutured manifold associated to
$(S,P,h)$ is isotopic to $(M,\G)$;  and that $\xi$ is compatible with
$(S,P,h)$ --- all in the sense that we defined in this paper. In the
rest of this section we prove these claims and Lemma~\ref{etnyre}
to obtain a proof of Theorem~\ref{giroux}.

The tubular portion $T$ of $ -\p N(K)$ in Theorem~\ref{honkazmat}(A)(i) is
split by its dividing set into positive and negative
regions, with respect to the orientation of $\p (M \setminus
N(K))$. Let $P$ be the positive region. Note that the negative
region $ T \setminus P$ is diffeomorphic to $P$. Since $(M, \G)$
is assumed to be a (balanced) sutured manifold, $\p M$ is divided
into $R_+(\G)$ and $R_-(\G)$ by the suture $\G$. Let $R_+ =
R_+(\G) \setminus \cup_i D_i$, where $D_i$'s are
the components of $\p N(K) \cap \p M$
and let $S$ be the surface which is
obtained from $\overline{R_+ }$ by attaching the positive region
$P$. If we denote the dividing set of $T$ by $A = \p P\cap \p S$,
then it is easy to see that
$$N(K) \cong (P \times [0,1])/\sim $$ where $(x,t) \sim (x,t')$
for $x \in A$ and $t, t' \in [0,1]$, such that the dividing set of
$\p N(K) $ is given by $ \p P \times \{1/2\}$.

In \cite{hkm1}, Honda, Kazez and Mati\'{c} observed that
$$\overline{M \setminus N(K)} \cong (S \times [-1,0]) / \sim $$
where $(x,t) \sim (x,t')$ for $x \in \p S$ and $t, t' \in [-1,0]$,
such that the dividing set of $\overline{M \setminus N(K)} $ is
given by $ \p S \times \{0\}$.

Moreover the embedding $h: P \to S$ which is obtained by first
pushing $P$ across $N(K)$ to $ T \setminus P \subset \p (M
\setminus N(K))$, and then following it with the identification of
$\overline{M \setminus N(K)}$ with $(S \times [-1,0])/ \sim$ is
called the monodromy map in the Honda-Kazez-Mati\'{c} description
of a partial open book decomposition.

In conclusion, we see that the triple $(S,P,h)$ satisfies the
conditions in Definition~\ref{Pob}:

$(1)$ The compact oriented surface $S$ is connected since we assumed
that $M$ is connected and it is clear that $\p S \neq \emptyset$.

$(2)$ The surface $P$ is a proper subsurface of $S$ such that $S$
is obtained from $\overline{S \setminus P}$ by successively
attaching $1$-handles by construction.

$(3)$ The monodromy map $h:P \to S$ is an embedding such that $h$
fixes $A = \p P \cap \p S$ pointwise.

Next we observe that $N(K)$  (resp. $\overline{M \setminus N(K)}$)
corresponds to $N$ (resp. $H$) in our construction of the balanced
sutured manifold associated to a partial open book decomposition
proceeding Definition~\ref{Pob}. The monodromy map $h$ amounts to
describing how $N=N(K)$ and $H =\overline{M \setminus N(K)}$ are
glued together along the appropriate subsurface of their
boundaries. This proves that the balanced sutured manifold
associated to $(S,P,h)$ is diffeomorphic to $(M,\G)$.

\begin{Lem}\label{compatible} The contact structure $\xi$ in
Theorem~\ref{honkazmat} is compatible with the \pob $(S,P,h)$
described above.
\end{Lem}

\begin{proof} We have to show that the contact structure $\xi$
in Theorem~\ref{honkazmat} satisfies the conditions $(1),(2)$ and
$(3)$ stated in Proposition~\ref{torisu}  with  respect to the  \pob
$(S,P,h)$ described above. We already observed that $N= N(K) $ and
$H= \overline{M \setminus N(K)}$. Then

$(1)$ The restrictions
of the contact structure $\xi$ onto $N(K)$ and $\overline{M
\setminus N(K)}$ are tight by conditions (A)(iii) and (B)(ii) of
Theorem~\ref{honkazmat}, respectively. This is because in either
case one obtains a standard contact $3$-ball or a disjoint union
of standard contact $3$-balls by cutting the manifold along a
collection of compressing disks each of whose boundary
geometrically intersects the dividing set exactly twice, and hence the dividing set of each of these compressing disks is a single boundary-parallel arc (see \cite[Corollary 2.6 (2)]{h}).

$(2)$
$\p H = \p (M\setminus N(K)) = (\p M \setminus \cup_i D_i) \cup
T$\ is convex by the convexity of $\p M$ and the convexity of $T$
(condition (A)(i) in Theorem~\ref{honkazmat}). Its dividing set is
the union of those of $\p M \setminus \cup_i D_i$ and $T$, hence
it is isotopic to $(\overline{\p S \setminus \p P}) \times \{ 0\}
\cup A \times \{ 0\} = \p S \times \{ 0\}$.

$(3)$ $\p N = \p
N(K) = \cup_i D_i \cup T$ is convex by the convexity of $D_i
\subset \p M$ and the convexity of $T$. Its dividing set is the
union of those of $D_i$'s and $T$, hence it is isotopic to
$(\overline{\p P \setminus \p S}) \times \{ 1/2\} \cup A \times \{
0\} = \p P \times \{ 1/2\}$.
\end{proof}

The following lemma is the only remaining ingredient in the proof
of Theorem~\ref{giroux}.

\begin{Lem}\label{etnyre}
Let $(S,P,h)$ be a partial open book decomposition, $(M,\G)$ be
the balanced sutured manifold associated to it, and $\xi$ be a
compatible contact structure. Then
$(S,P,h)$ is
given by the Honda-Kazez-Mati\'c description.
\end{Lem}
\begin{proof}

Consider the graph $K$ in $P$ that is obtained by gluing the core
of each $1$-handle in $P$ (see Figure~\ref{handles} for example).

\begin{figure}[ht]
 \relabelbox \small {
 \centerline{\epsfbox{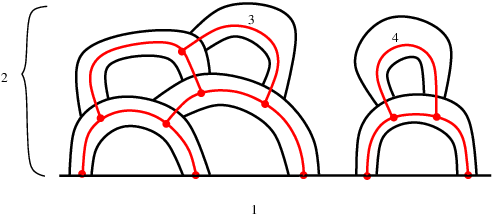}}}

 \relabel{1}{$ S \setminus P $}
 \relabel{2}{$P$}
 \relabel{3}{$K$}
 \relabel{4}{$K$}

\endrelabelbox
       \caption{ Legendrian graph $K$ in $P$}
       \label{handles}
\end{figure}

It is clear that $P$ retracts onto $K$. We will denote $K \times
\{1/2\} \subset P \times \{ 1/2\}$ also by $K$. We can first make
$P \times \{1/2\}$ convex and then Legendrian realize $K$ with
respect to the compatible contact structure $\xi$ on $N \subset
M$. This is because each component of the complement of $K$ in $P$
contains a boundary component (see \cite[Remark 4.30]{e}). Hence
$K$ is a Legendrian graph in $(M,\xi)$ with endpoints in $\p P
\times \{1/2 \} \setminus \p S \times \{ 0 \} \subset \G \subset
\p M$ such that $N=P \times [0,1] / \sim $ is a neighborhood
$N(K)$ of $K$ in $M$. Then all the conditions except (A)(i) in
Theorem~\ref{honkazmat} on $N(K) =N$ and $\overline{M\setminus
N(K)} = H$
are satisfied because of the way we constructed $\xi$ in
Proposition~\ref{torisu}. Since $\p N$ is convex
$T$ is also convex.
It remains to check that the boundary of the tubular portion $T$ of
$N$ is Legendrian. Note that each component of this boundary $\p
D_i = \p (\g_i \times [0,1]) \subset \p N$ is identified with
$\gamma_i = \g_i \times \{0\} \cup h(\g_i) \times \{ -1\}$ in the
convex surface $\p H = S\times \{0\} \cup -S \times \{ -1\}$.
Since each $\gamma_i$ intersects the dividing set $\G_{\p H} = S
\times \{ 0\}$ of $\p H$ transversely at two points $\p c_i \times
\{ 0\}$, the set $\{\gamma_1, \gamma_2, \ldots, \gamma_n\}$ is
non-isolating in $\p H$ and hence we can use the Legendrian
Realization Principle to make each $\gamma_i$ Legendrian.
\end{proof}

\noindent {\emph{Proof of Theorem~\ref{giroux}.} By
Proposition~\ref{torisu} each \pob is compatible with a unique
compact contact $3$-manifold with convex boundary up to contact
isotopy. This gives a map from the set of all partial open book
decompositions to the set of all compact contact $3$-manifolds
with convex boundary and by Remark~\ref{unique} this map descends
to a map from the set of isomorphism classes of all partial open
book decompositions  to the set of isomorphism classes of all
compact contact $3$-manifolds with convex boundary. Moreover by
Lemma~\ref{invariance} this gives a well-defined map $\Psi$ from
the isomorphism classes of all partial open book decompositions
modulo positive stabilization to that of isomorphism classes of
compact contact $3$-manifolds with convex boundary. On the other
hand, Honda-Kazez-Mati\'c description gives a well-defined map
$\Phi$ in the reverse direction by \cite[Theorems~1.1 and 1.2]
{hkm1}. Furthermore, $\Psi \circ \Phi$ is identity by
Lemma~\ref{compatible} and $\Phi \circ \Psi$ is identity by
Lemma~\ref{etnyre}. \hfill $\Box$

\section{Examples}

Below we provide examples of abstract partial open books which correspond to some basic contact $3$-manifolds with boundary. These examples were previously appeared in \cite{eo} where their contact invariants were calculated.

\v

\begin{exa}\label{ex1} Let $S$ be an annulus, $P$ be a regular
neighborhood of $r$ disjoint arcs
connecting the two distinct boundary components of $S$ as in
Figure~\ref{ex1-1}, and the
monodromy $h$ be the inclusion of $P$ into $S$. The partial open book $(S,P,h)$ is compatible with the contact structure
obtained by removing $r$ disjoint standard contact open $3$-balls
from the unique (up to isotopy) tight contact structure
$\xi_{std}$ on $S^1 \times S^2$.

\begin{figure}[ht]
  \relabelbox \small {
  \centerline{\epsfbox{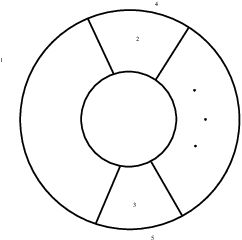}}}
  \relabel{1}{{$S$}}
  \relabel{2}{{$P_1$}}
 \relabel{3}{{$P_r$}}
 
  \endrelabelbox
        \caption{The annulus $S$, $r$ components $P_1, \dots , P_r$
        of $P$ 
                 in Example~\ref{ex1}.}
        \label{ex1-1}
\end{figure}
\end{exa}

\begin{exa}[\textbf{Standard contact $3$-ball}]\label{ex2} Let $S$ and
$P$ be as in Example~\ref{ex1} for $r=1$, and the monodromy $h$ be
the restriction (to $P$) of a {\emph{right}}-handed Dehn twist
along the core of $S$.
The contact $3$-manifold $(M, \G, \xi)$ compatible with this partial open book is
the standard contact $3$-ball. Here
the Legendrian graph $K$ which satisfies the conditions in
Theorem~\ref{honkazmat} is a single arc in $B^3$ connecting two
distinct points on $\G$ as depicted in Figure~\ref{poball}. The
complement $H$ of a regular neighborhood $N = N(K)$ in the
standard contact $3$-ball $B^3$ is a solid torus with two parallel
 dividing curves (see Figure~\ref{divcurve}) on $\p H$ which are homotopically
 nontrivial inside $H$. Here a meridional disk in $H$ will serve as the
 required compressing disk  $D_1^\d$ for $H$ in Theorem~\ref{honkazmat}
 $(B)$. On the other hand, $N$ is already a standard contact $3$-ball.
This shows in particular that the standard contact $3$-ball can be
obtained from a tight solid torus $H$ by attaching a tight
$2$-handle $N$. 

\begin{figure}[ht]
  \relabelbox \small {
  \centerline{\epsfbox{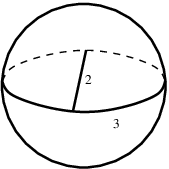}}}

  \relabel{2}{{$K$}}
 \relabel{3}{{$\G$}}

  \endrelabelbox
        \caption{The Legendrian arc $K$ in the standard contact $3$-ball.}
        \label{poball}
\end{figure}

\begin{figure}[ht]
  \relabelbox \small {
  \centerline{\epsfbox{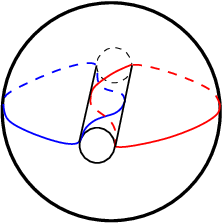}}}

  \endrelabelbox
        \caption{The dividing curves  on $\p H$.}
        \label{divcurve}
\end{figure}
\end{exa}

\begin{exa}[\textbf{Standard neighborhood of an overtwisted
disk}]\label{ex3} Let $(S,P,h)$ be the partial open book decomposition shown
in Figure~\ref{ex4-1}. This is the partial open book
considered in \cite[Example 1]{hkm1} which is compatible with
the standard neighborhood of an overtwisted disk.

Here we observe that by Proposition~\ref{torisu}, $(M, \G, \xi)$
is obtained by gluing a pair of compact connected contact $3$-manifolds
with convex boundaries, namely $(H, \G_{\p H} , \xi\vert_H)$ and
$(N, \G_{\p N} , \xi\vert_N)$, along parts of their boundaries.
We know that
$$H=(S \times [-1,0])/\sim $$ where $S$ is an annulus and $(x,t)
\sim (x,t')$ for $x \in \p S$ and $t, t' \in [-1,0]$. There is a
unique (up to isotopy) compatible tight contact structure on $H$
whose dividing set $\G_{\p H}$ on $\p H$  is $\p S \times \{0\}$
(cf. Proposition~\ref{torisu}). Hence $(H, \G_{\p H}, \xi\vert_H)$
is a solid torus carrying a tight contact structure where $\G_{\p
H}$ consists of two parallel curves on $\p H$ which are
homotopically nontrivial in $H$. We observe that when we cut $H$
along a compressing disk we get a standard contact $3$-ball $B^3$
with its connected dividing set $\G_{\p B^3}$ on its convex
boundary. Note that $\G_{\p B^3}$ is obtained by ``gluing" $\G_{\p
H}$ and the dividing set on the compressing disk. 
Similarly we know that $N=(P \times [0,1])/\sim$, where $(x,t)
\sim (x,t')$ for $x \in A$ and $t, t' \in [0,1]$. There is a
unique (up to isotopy) compatible tight contact structure on $N$
whose dividing set $\G_{\p N}$ on $\p N$  is $\p P \times \{1/2\}$
(cf. Proposition~\ref{torisu}). We observe that $(N, \G_{\p N},
\xi\vert_N)$ is the standard contact $3$-ball.
\end{exa}

\bibliographystyle{amsalpha}

\end{document}